\documentclass[11pt]{amsart}

\RequirePackage[l2tabu,orthodox]{nag}
\usepackage[T1]{fontenc}
\usepackage[utf8]{inputenc}
\usepackage{lmodern}
\usepackage[english]{babel}
\usepackage{color} 
\usepackage[usenames,dvipsnames]{xcolor}
\usepackage{amsrefs}
\usepackage{amsthm,amssymb,amsmath}
\usepackage[pdftex]{graphicx}
\usepackage{enumerate}
\usepackage[inline]{enumitem}
\usepackage{url}
\usepackage{faktor}

\usepackage{multicol}
\usepackage{mathrsfs}

\numberwithin{equation}{section}

\let\oldtocsection=\tocsection

\let\oldtocsubsection=\tocsubsection 

\let\oldtocsubsubsection=\tocsubsubsection

\renewcommand{\tocsection}[2]{\hspace{0em}\oldtocsection{#1}{#2}}
\renewcommand{\tocsubsection}[2]{\hspace{2em}\oldtocsubsection{#1}{#2}}
\renewcommand{\tocsubsubsection}[2]{\hspace{4em}\oldtocsubsubsection{#1}{#2}}

\newtheorem{thm}{Theorem}[section]

\newtheorem{ques}[thm]{Question}

\newtheorem{lem}[thm]{Lemma}

\newtheorem{cor}[thm]{Corollary}
\newtheorem{intro-cor}[introthm]{Corollary}

\newtheorem{conj}[thm]{Conjecture}
\theoremstyle{definition}
\newtheorem{dfn}[thm]{Definition}

\newtheorem*{main question}{Main Question}

\theoremstyle{remark}
\newtheorem{rem}[thm]{Remark}
\newtheorem{ex}[thm]{Example}

\usepackage[margin=2.5cm]{geometry}

\usepackage{indentfirst}
\usepackage{microtype}

\usepackage[colorlinks=true,linkcolor=blue,citecolor=magenta]{hyperref}

\usepackage{blindtext}

\usepackage{nameref}

\DeclareMathOperator{\homeo}{\mathsf{Homeo}}
\DeclareMathOperator{\Diff}{\mathsf{Diff}}

\DeclareMathOperator{\PL}{\mathsf{PL}}

\newcommand{\cN}{\mathcal{N}}

\newcommand{\R}{\mathbb{R}}
\newcommand{\Q}{\mathbb Q}
\newcommand{\N}{\mathbb N}
\newcommand{\Z}{\mathbb Z}
\newcommand{\T}{\mathbb S^1}

\renewcommand{\setminus}{\smallsetminus}

\DeclareMathOperator{\Int}{\mathsf{Int}}

\newcommand{\rot}{\mathsf{rot}}

\title[Non-smoothability for a class of groups of PL homeomorphisms of the interval]{Non-smoothability for a class of groups of piecewise linear homeomorphisms of the interval}

\author{Michele Triestino}
\date{\today}

\begin{document}
\begin{abstract}

For a certain class of groups of piecewise linear homeomorphisms of the interval, we prove that they admit no sufficiently regular faithful action on the line. Building on previous work of Brum, Matte Bon, Rivas, and the author [arXiv:2104.14678], the new ingredient is an observation from a recent work of Hyde and Tatch Moore [arXiv:2103.14911], which allows to reduce the problem to the case of the circle, and then apply Herman--Yoccoz theory.

	\smallskip
	
	{\noindent\footnotesize \textbf{MSC\textup{2020}:} Primary 37C85, 57M60. Secondary 37E05, 37E10.}

\end{abstract}

\maketitle

\section{Introduction} The class of groups of piecewise linear (PL) homeomorphisms\footnote{We will always tacitly consider only \emph{orientation-preserving} homeomorphisms.} of the real line is a prolific source of examples of groups with striking properties. The emblematic example is Thompson's group $F$, a group admitting a presentation with two generators and two relations, whose derived subgroup is simple. The group $F$ consists of all PL homeomorphisms of $[0,1]$ that satisfy some arithmetic conditions. More precisely, for a given real $\alpha\in \R$, we let $F_\alpha$ be the group of all PL homeomorphisms of $[0,1]$, which are locally of the form $x\mapsto \alpha^kx+b$, for some $k\in \Z$ and $b\in \Z[\alpha,\alpha^{-1}]$, and whose breaks of derivatives are also in $\Z[\alpha,\alpha^{-1}]$. With this notation, we have $F=F_2$. 

We will be interested in regularity properties of actions of groups of PL homeomorphisms. It was first observed by Ghys and Sergiescu \cite{GhysSergiescu} that the natural action of $F$ on $[0,1]$ can be conjugated to an action by $C^\infty$ diffeomorphisms (this can be now easily deduced from the work of Kim, Koberda, and Lodha \cite{KKL}, using Brin's ``2-chain lemma''). The method of Ghys and Sergiescu can be extended to the groups $F_n$, where $n\ge 3$ is an integer, but actually $F_n$ embeds in $F$ (see for instance Brin and Guzm\'an \cite{BrinGuzman}).
On the other hand, groups defined by higher rank slope groups, such as the Brown--Thompson--Stein groups $F_{n_1,\ldots,n_k}$ with $k\ge 2$, admit no faithful action on $[0,1]$ of class $C^r$ for any $r>1$, as  proved by Brum, Matte Bon, Rivas, and the author in \cite{brum2021locally} (the case $r=2$ was previously discussed in the work by Bonatti, Lodha, and the author \cite{BLT}, and it applies to a larger class of groups). This implies for instance that Brown--Thompson--Stein groups do not embed in $F$, a result that was first obtained by Lodha \cite{Coherent}. Conversely, the results from \cite{BLT,brum2021locally} are partially motivated by the following conjecture.

\begin{conj}\label{conj.main}
	Let $G\le \PL([0,1])$ be a finitely generated group of PL homeomorphisms of $[0,1]$, whose standard action is conjugate to a $C^\infty$ action. Then $G$ embeds in Thompson's group $F$.
\end{conj}

No counterexample is known, even if one replaces the $C^\infty$ assumption by $C^1$, or conjugacy by semi-conjugacy. Surprisingly, \emph{rational} (non integer) slope Thompson groups have not been considered so far in the literature on Thompson groups, and in such case both statements in Conjecture \ref{conj.main} are unsolved.

\begin{ques}
	Take a rational $r\in \Q\setminus \Z$, $r>1$.
	\begin{enumerate}
		\item Does the group $F_r$ embed in $F$?
		\item Is the standard action of $F_r$ conjugate to a $C^\infty$ action?
	\end{enumerate}
\end{ques}

Another particular, but relevant case for which the available techniques do not apply is given by \emph{irrational slope Thompson groups}, namely the groups $G=F_\alpha$ for irrational $\alpha$. It was recently remarked by Hyde and Tatch Moore \cite{hyde2021subgroups} that such groups do not embed in $F$, introducing an elementary technique that inspired our work. 
More precisely, the aim of this note is to use Herman--Yoccoz theory to show that some groups of homeomorphisms of the line cannot be semi-conjugate to groups which are of class $C^r$, for some appropriate $r>2$. The main technical result will be Theorem \ref{t.keythm} below. In relation to Conjecture \ref{conj.main}, we are able to obtain the following.

\begin{cor}\label{cor.Falpha}
	Let $\alpha\in \R\setminus \Q$ be an irrational number of bounded type (for instance, a quadratic irrational number). Then the group $F_\alpha$ admits no faithful action on $[0,1]$ of class $C^r$, for any $r>2$. 
	
	More generally, if $\alpha\in \R$ satisfies a Diophantine condition $\mathcal D_\delta$ (see the definition below), then the group $F_\alpha$ admits no faithful action on $[0,1]$ of class $C^r$, for any $r>2+\delta$. 
\end{cor}

Results of similar spirit were obtained by Liousse in \cite{liousse}, were she was considering groups of PL homeomorphisms of the circle containing elements of irrational rotation number.

\section{Herman--Yoccoz theory}
Throughout the paper, given a one-manifold $X$, we will denote by $\homeo_0(X)$ the group of orientation-preserving homeomorphisms of $X$, and by $\Diff_0^r(X)$ the group of orientation-preserving $C^r$ diffeomorphisms of $X$.

One of the most fundamental results in one-dimensional dynamics is the celebrated Denjoy theorem, which states that any circle diffeomorphism of class $C^2$ without periodic orbits is minimal, namely every orbit is dense. By Poincaré theory, this implies that any such diffeomorphism $f$ is topologically conjugate to an irrational rotation of the circle, whose angle is determined by the \emph{rotation number} of $f$, defined by
\[
\rot (f)=\lim_{n\to\infty}\frac{F^n(0)}{n}\quad\pmod{\Z},
\]
where $F:\R\to\R$ is any lift of $f$ to a diffeomorphism of the real line. The $C^2$ regularity assumption is far from optimal: for instance, the conclusion of Denjoy theorem also holds for homeomorphisms of the circle which belong to Herman's \emph{class $\mathrm{P}$} \cite[\S VI.4]{Herman}, formed by all homeomorphisms of the circle $g$ whose left derivative $D_-g$ exists everywhere and such that $\log D_-g$ is a function of bounded variation\footnote{This is not Herman's definition of class $\mathrm P$, but the condition is equivalent, as noted by Herman in the remark right after the definition in \cite[\S VI.4]{Herman}.}. To fix notation, we will also denote by $D_+g$ the right derivative of a map $g$, and by $\sigma(g)=D_+g/D_-g$ its \emph{jump} of derivatives.

\begin{thm}[Herman]\label{t.herman}
	Let $f\in \homeo_0(\T)$ be a circle homeomorphism in the class $\mathrm{P}$ and of irrational rotation number $\alpha=\rot(f)$. Then $f$ is $C^0$ conjugate to the rotation $R_\alpha$. 
\end{thm}

Here we will be interested in homeomorphisms $g$ which are piecewise of class $C^r$, for some $r>2$, and whose left derivative is bounded away from $0$. Such maps belong to the class $\mathrm{P}$, and include all PL homeomorphisms. The scope of \emph{Herman--Yoccoz theory} is to go beyond Denjoy theorem and determine which conditions on the regularity of $f$ and the arithmetic nature of the rotation number $\alpha=\rot(f)$ guarantee a better regularity for the maps which conjugates $f$ to $R_\alpha$.
For this, recall that a real number $\alpha$ satisfies a \emph{Diophantine condition} $\mathcal D_\delta$ if there exists a constant $C>0$ such that $|q\alpha-p|>Cq^{-1-\delta}$ for every $p/q\in\mathbb Q$.

\begin{ex}When $\alpha$ is of \emph{bounded type}, namely if the terms appearing in the continued fraction expansion of $\alpha$ are uniformly bounded (this happens for instance for quadratic irrational numbers, for which the continued fraction expansion is periodic), then $\alpha$ satisfies the Diophantine condition $\mathcal D_\delta$ for any $\delta>0$.
\end{ex}

\begin{rem}\label{r.SL2}
	The Diophantine condition $\mathcal D_\delta$ is preserved by the action of $\mathsf{PSL}(2,\Z)$ by fractional linear transformations, as such operations do not change the tail of the continued fraction expansions.
\end{rem}

Among the sharpest statements in Herman--Yoccoz theory, we have the following result by Khanin and Teplinsky \cite{KhaninTeplisnky}.

\begin{thm}[Khanin--Teplinsky]\label{t.KhaninTeplinsky}
	For $r>2$, let $f \in \Diff_0^r(\T)$ be a diffeomorphism whose rotation number satisfies the Diophantine condition $\mathcal D_\delta$. If $r>2+\delta$, then $f$ is $C^{r-\delta-1}$ conjugate to the corresponding irrational rotation. 
\end{thm}

\section{Local irrational rotations}

In their work \cite{hyde2021subgroups}, Hyde and Tatch Moore introduced a (semi-)conjugacy invariant for group actions on the real line, using rotation number, which is the most fundamental (semi-)conjugacy invariant for circle homeomorphisms. Before introducing this, given non-empty open real intervals $X=(a,b)$ and $X'=(a',b')$, one says that two actions $\rho:G\to \homeo_0(X)$ and $\rho':G\to \homeo_0(X')$ are \emph{(positively) semi-conjugate}, if there exists a non-decreasing map $\varphi:X \to X'$, which is proper (that is, $\lim_{x\to a} \varphi(x)=a'$ and similarly $\lim_{x\to b} \varphi(x)=b'$) and such that $\varphi\rho(g)=\rho'(g)\varphi$ for any $g\in G$. We call the map $\varphi$ a \emph{semi-conjugacy} between $\rho$ and $\rho'$.
When the semi-conjugacy is an orientation-preserving homeomorphism, one says that $\rho$ and $\rho'$ are \emph{(positively) conjugate}. With this definition, semi-conjugacy is an equivalence relation on the space of actions $\rho:G\to \homeo_0(\R)$ of $G$ on the line. Given two circles $C$ and $C'$, we say that two actions $\rho:G\to \homeo_0(C)$ and $\rho':G\to \homeo_0(C')$ of a group $G$ are (semi-)conjugate, if one can lift $\rho$ and $\rho'$ to (semi-)conjugate actions on their universal covers of the associated cyclic central extension of $G$. In such case, a map $\widetilde\varphi$ realizing the semi-conjugacy is equivariant with respect to the monodromy actions of $\Z$, and thus defines a map $\varphi:C\to C'$ such that $\varphi \rho(g)=\rho'(g)\varphi$ for every $g\in G$. We call the map $\varphi$ a \emph{semi-conjugacy} between $\rho$ and $\rho'$.

Assume that $X$ and $X'$ are one-manifolds, $\rho:G\to \homeo_0(X)$ and $\rho':G\to \homeo_0(X')$  are semi-conjugate actions of a group $G$, with the action $\rho'$ being minimal. Let $\varphi$ be a semi-conjugacy for the actions. Then, minimality of $\rho'$ ensures that $\varphi$ is surjective, and thus continuous. It follows that the subset
\[U=\{x'\in X':\varphi^{-1}(x')\text{ is a singleton}\}\]
is $\rho'(G)$-invariant and its complement is at most countable. We will call $U$ the \emph{set of bi-univocality} of the semi-conjugacy $\varphi$.

A result of Ghys \cite{Ghys-bounded} says that the semi-conjugacy class of an action of a group $G$ on the circle is completely determined by the so-called \emph{bounded Euler class}, which can be computed in terms of rotation numbers of the elements (more precisely, in terms of \emph{translation numbers}). In contrast, no reasonable complete invariant is available for semi-conjugacy classes of actions on the real line (see the discussion in \cite{brum2021locally}), but the invariant defined by Hyde and Tatch Moore seems to be quite powerful, for some specific classes of groups.
Let us explain the general idea. Given a subgroup $G\le \homeo_0(\R)$, one can consider all the piecewise continuous maps $T:I\to I$ which are locally defined by elements of $G$, and where $I\subset \R$ is an arbitrary closed interval. More precisely, we require that there exists a decomposition $I=I_1\cup \cdots \cup I_n$ of the interval into $n$ adjacent subintervals with pairwise disjoint interior, and elements $g_1,\ldots,g_n\in G$ such that $g_i(I_i)\subset I$ for any $i\in \{1,\ldots,n\}$, and such that $T(x)=g_i(x)$ for any $x\in \Int(I_i)$. The basic, but fundamental observation, is that any (semi-)conjugacy of the action of $G$ induces a piecewise continuous map $\widetilde{T}$ which is (semi-)conjugate to $T$.

\begin{rem}Actually, Hyde and Tatch Moore were only considering maps $T$ which induce a circle homeomorphism (after identifying the endpoints of $I$), but it can be interesting, although we have not yet found significant applications, to look at more general locally defined transformations $T$, such as generalized interval exchange maps of $I$, or doubling maps of the circle obtained by identifying the endpoints of $I$. For instance, results from \cite{HmiliLiousse,CentralizersVn} allow to describe the dynamics of affine interval exchange transformations locally defined by elements of Thompson's group $F$, and for certain generalized interval exchange maps some Herman--Yoccoz theory has been developed by Ghazouani and Ulcigrai \cite{ghazouani2021priori}.\end{rem}

Here, as in \cite{hyde2021subgroups}, we will be only interested in the case of local irrational rotations.

\begin{dfn}\label{d.locrot}For $X=(a,b)$, let $G\le \homeo_0(X)$ be a subgroup. A \emph{local rotation of angle $\alpha$} in the group $G$ consists of a point $y\in X$ and a pair of elements $f,g\in G$ with the following properties:
\begin{enumerate}
	\item $g(y)=:x<y<z:=f(y)$;
	\item $gf(y)=fg(y)$;
	\item\label{i.epsilon} there exists $\varepsilon>0$ such that $Df(t)=1$ for every $t\in (x-\varepsilon,y+\varepsilon)$ and $Dg(t)=1$ for every $t\in (y-\varepsilon,z+\varepsilon)$;
	\item $\dfrac{fg(y)-x}{z-x}=\alpha\pmod\Z$.
\end{enumerate}
We denote by $(y;f,g)$ the local rotation.
\end{dfn}

\begin{rem}\label{r.locrot}
	In other terms, the pair $(f,g)$ defines an interval exchange transformation $T$ of the interval $[x,z]$ which is affinely conjugate to the rotation of angle $\alpha$ on $\R/\Z$ (although the third condition is slightly stronger).\end{rem}

This is our main technical result.

\begin{thm}\label{t.keythm}
	For $X=(a,b)$, let $G\le \homeo_0(X)$ be a subgroup with the following properties:
	\begin{enumerate}
		\item the action of $G$ on $X$ is minimal;
		\item $G$ contains a local rotation $(y;f,g)$ of Diophantine angle $\alpha\in \mathcal D_\delta$;
		\item $G$ contains an element $h$ with two fixed points in $(x,z)$, at which it is not differentiable.
	\end{enumerate}
	Then, for any $r>2+\delta$, the standard action of $G$ is not semi-conjugate to any $C^{r}$ action on $\R$.
\end{thm}

\begin{proof}
	Assume by contradiction there exists $\widetilde \rho:G\to \Diff^r_0(\R)$ which is semi-conjugate to the standard action $\rho_0$ of $G$.
	We first claim that $\widetilde \rho$ has to be \emph{conjugate} to $\rho_0$. Indeed, as we are assuming that the action of $G$ is minimal, any semi-conjugacy must be realized by a \emph{continuous} non-decreasing map $\varphi:\R\to X$, such that $\varphi \widetilde \rho=\rho_0\varphi$.
	Write $I=[x,z]$ and consider the interval $\widetilde I:=\varphi^{-1}(I)$. We also write $\widetilde f=\widetilde \rho(f)$ and $\widetilde{g} =\widetilde \rho(g)$. Because of condition \eqref{i.epsilon} for a local rotation, up to replace the point $y$ by a nearby point $y'\in (y-\varepsilon,y+\varepsilon)$, we can assume that $\varphi^{-1}(y)$ is a singleton, namely, we can assume that $y$ belongs to  $U$, the set of bi-univocality of $\varphi$ (which is dense). Moreover, as $U$ is $\rho_0(G)$ invariant, this ensures that the points $\widetilde x=\varphi^{-1}(x)$ and $\widetilde z=\varphi^{-1}(z)$ are well-defined, and we have $\widetilde I=[\widetilde x,\widetilde{z}]$. Write $\widetilde y=\varphi^{-1}(y)$, then the map $\widetilde T:\widetilde I\to \widetilde I$, defined by
	\[
	\widetilde T(t)=\left\{\begin{array}{lr}
	\widetilde f(t)& \text{if }t\in [\widetilde x,\widetilde y],\\[.5em]
	\widetilde g(t)& \text{if }t\in (\widetilde y,\widetilde z],
	\end{array} \right.
	\]
	induces a homeomorphism of the circle obtained by identifying the endpoints of $\widetilde I$, and it is semi-conjugate (by the composition of $\varphi$ with the affine map coming from Remark \ref{r.locrot}) to the rotation of angle $\alpha$. Keep denoting such homeomorphism by $\widetilde T$, and remark that $\widetilde T$ is piecewise of class $C^r$, with possible breaks at $\widetilde x\sim\widetilde z$ and $\widetilde y$, at which $D_-\widetilde T$ is bounded away from $0$. After Denjoy theorem for the class  $\mathrm{P}$, the semi-conjugacy induced by $\varphi$ must be a conjugacy. By minimality of $\rho_0$, we deduce that the whole map $\varphi:\R\to X$ must define a conjugacy between $\widetilde \rho$ and $\rho_0$.
	
	We want to apply Herman--Yoccoz theory to deduce that $\varphi$ must be of class $C^{r-\delta-1}$. However, the possible two breaks of $\widetilde T$ do not allow us to use Theorem \ref{t.KhaninTeplinsky} directly. Note that the breakpoints $\widetilde x$ and $\widetilde y$ satisfy $\widetilde T(\widetilde{y})=\widetilde x$, so they belong to the same $\widetilde T$-orbit. Circle maps whose breakpoints belong to the same orbit have been studied by several authors. For the interested reader, we suggest the relatively recent work of Khanin and Koci\'c in \cite{KhaninKocic}. For our purposes, we will use a result by Adouani and Marzougui \cite[Theorem 2.1]{MR2427509}, extending previous work of Dzhalilov and Liousse \cite{Dzhalilov_2006}. We detail this.
	
	The piecewise $C^r$ circle homeomorphism $\widetilde{T}$ has $\{\widetilde{y},\widetilde T(\widetilde y)\}$ as only possible breakpoints. In such elementary case, the definition of condition $D_r$ by Adouani and Marzougui \cite[Définition 1.7]{MR2427509}, means that  $\widetilde T^2$ is of class $C^{\lfloor r\rfloor}$ at $\widetilde y$. 
	To check this, observe that $\widetilde T^2= \widetilde f\widetilde g$ on some small right neighborhood of $\widetilde y$, and $\widetilde T^2= \widetilde g\widetilde f$ on some small left neighborhood of $\widetilde y$. By condition \eqref{i.epsilon} in the definition of a local rotation, these compositions coincide on some small neighborhood of $\widetilde{y}$: we have $fg=gf$ on a neighborhood of $y$, so after conjugation, the same holds for $\widetilde f$ and $\widetilde g$ on a neighborhood of $\widetilde y$. As the maps $\widetilde f$ and $\widetilde g$ are $C^r$ at $\widetilde y$, so is $\widetilde T^2$.
	
	Using \cite[Theorem 2.1]{MR2427509} (and its proof, see also the proof of \cite[Corollary 1.2]{KhaninKocic}), we get that $\widetilde T$ is  conjugate to a $C^r$ circle diffeomorphism by a homeomorphism which is piecewise $C^r$ and has one break. In particular, after Theorem \ref{t.KhaninTeplinsky}, it is conjugate to the corresponding irrational rotation by a map which is piecewise $C^1$ and has one break 
	(see Lemma \ref{l.reg_conjugation} for a direct proof in the $C^{2+\tau}$ case). Note that $\varphi^{-1}$ and $\psi$ differ of a rotation in restriction  to $[x,z]$, so $\varphi^{-1}$ has the same regularity as $\psi$.
	Now, taking the element $h$ from the statement, we deduce that the conjugate $\varphi^{-1} h\varphi$ has a discontinuity point for the derivative on the interval $(\widetilde x,\widetilde z)$.
	\qedhere 
\end{proof}

\section{Irrational slope Thompson groups}

Here we apply Theorem \ref{t.keythm} to show that the irrational slope Thompson groups $F_\alpha$ with Diophantine $\alpha\in \mathcal D_\delta$, do not admit $C^r$ faithful actions, for any sufficiently large $r>2$ (Corollary \ref{cor.Falpha}).

We first show that $F_\alpha$ contains local rotations of Diophantine angle. Our original construction was erroneous, and we thank the referee for fixing it. We first need a basic result from the monograph of Bieri and Strebel \cite[Theorem A4.1]{BieriStrebel}. To state it, given a multiplicative subgroup $\Lambda\le (0,+\infty)$, we write $A=\Z[\Lambda]$ for its group ring, and 
\[I\Lambda\cdot A:=\langle(\lambda-1)a: \lambda\in \Lambda, a\in A\rangle,\]
seen as $\Z[\Lambda]$-submodule of $A$. We denote by $G(\R;A,\Lambda)$ the \emph{Bieri--Strebel group} of all PL homeomorphisms of $\R$ with slopes in $\Lambda$ and breakpoints in $A$.

\begin{thm}[Bieri--Strebel]\label{t.BScriterion}
	Let $\Lambda\le (0,+\infty)$ be a multiplicative subgroup, and set $A=\Z[\Lambda]$. Let $a,a',c$ and $c'$ be points in $A$ with $a<c$ and $a'<c'$. Then, there exists an element $f\in G(\R;A,\Lambda)$ mapping $[a,c]$ onto $[a',c']$ if and only if
	\[
	c'-a'\equiv c-a\pmod{I\Lambda\cdot A}.
	\]
\end{thm}

Note that when $\Lambda=\langle \alpha\rangle$, we have $A=\Z[\alpha,\alpha^{-1}]$ and $I\Lambda\cdot A=(\alpha-1) A$. Note that when $\alpha=\frac{1+\sqrt{5}}{2}$ is the golden ratio, we have $I\Lambda\cdot A=A$, and thus there is no obstruction to find PL homeomorphisms as in the statement above. In general, there are obstructions, and our construction of local rotation will be slightly more involved than the one appearing in \cite[Section 7]{hyde2021subgroups} in the case of the golden ratio.

\begin{lem}
	For any $\alpha>1$ satisfying a Diophantine condition $\mathcal D_\delta$, the group $F_\alpha$ contains a local rotation of angle $\beta$ satisfying the same Diophantine condition $\mathcal D_\delta$.
\end{lem}

\begin{proof}
	We write $A=\Z[\alpha,\alpha^{-1}]$. To start with, take any point $y\in (0,1)\cap A$ and constants $\beta_1,\beta_2\in (\alpha-1)A$ such that 
	\[0<y-\beta_2<y<y+\beta_1<1.\]
	Note that this forces the conditions $0<\beta_1+\beta_2<1$ and $\beta_1,\beta_2>0$. Set $x:=y-\beta_2$ and $z:=y+\beta_1$, and take $\varepsilon>0$, with $\varepsilon\in A$, such that $0<x-\varepsilon$ and $z+\varepsilon<1$.
	Take the four points $a=0,c=x-\varepsilon$, and $a'=0,c'=x-\varepsilon+\beta_1$, and note that
	\[
	(c'-a')-(c-a)=c'-c=\beta_1\in (\alpha-1)A.
	\]
	Therefore, after Theorem \ref{t.BScriterion}, there exists a PL homeomorphism $f_1\in G(\R;\langle \alpha\rangle, A)$ sending $[0,x-\varepsilon]$ onto $[0,x-\varepsilon+\beta_1]$.
	Similarly, working with the four points $a=y+\varepsilon,c=1$ and $a'=y+\varepsilon+\beta_1,c'=1$, we can find a PL homeomorphism $f_2\in G(\R;\langle \alpha\rangle, A)$ sending $[y+\varepsilon,1]$ onto $[y+\varepsilon+\beta_1,1]$.
	Then, the map $f:[0,1]\to [0,1]$ defined by
	\[
	f(t)=\left\{
		\begin{array}{lr}
			f_1(t)&\text{if }t\in [0,x-\varepsilon],\\[.5em]
			t+\beta_1&\text{if }t\in (x-\varepsilon,y+\varepsilon),\\[.5em]
			f_2(t)&\text{if }t\in [y+\varepsilon,1]
		\end{array}
	\right.
	\]
	gives an element of $F_\alpha$ such that $f(t)=t+\beta_1$ for every $t\in (x-\varepsilon,y+\varepsilon)$. By an analogous argument, the condition $\beta_2\in (\alpha-1)A$ allows to find an element $g\in F_\alpha$ such that $g(t)=t-\beta_2$ for every $t\in (y-\varepsilon,z+\varepsilon)$.
	
	With such choices, we have
\[
\frac{fg(y)-x}{z-x}=\frac{\beta_1}{\beta_1+\beta_2},
\]
so we have found a local rotation $(y;f,g)$ of angle $\beta:=\frac{\beta_1}{\beta_1+\beta_2}$ in $F_\alpha$. We need to argue that we can choose $\beta_1$ and $\beta_2$ such that $\beta$ satisfies the same Diophantine condition $\mathcal{D}_\delta$ as $\alpha$. Note that if $\beta\in \mathcal D_\delta$, then $\beta^{-1}-1=\beta_2/\beta_1$ also satisfies the condition $\mathcal D_\delta$ (see Remark \ref{r.SL2}). As we are assuming $\alpha>1$, we can find $k\in \N$ such that
\[
0<y-(\alpha-1)\alpha^{-k}<y<y+(\alpha-1)\alpha^{-k-1}<1.
\]
Then the numbers $\beta_2=(\alpha-1)\alpha^{-k}$ and $\beta_1=(\alpha-1)\alpha^{-k-1}$ are in $(\alpha-1)A$, and thus satisfy all the requirements for our construction. Moreover, we have $\beta_2/\beta_1=\alpha\in \mathcal D_\delta$, so that $\beta=\frac1{1+\alpha}\in \mathcal{D}_\delta$, as desired.
\end{proof}

Now, after Theorem \ref{t.keythm}, we need to see that any $C^r$ faithful action on $[0,1]$ is semi-conjugate to the standard one. We proved this in our previous work in collaboration with Brum, Matte Bon, and Rivas \cite[Corollary 6.11]{brum2021locally}:

\begin{thm}
	For $X=(a,b)$, let $G\le \homeo_0(X)$ be locally moving with independent groups of germs at the endpoints. Then, for every faithful action $\rho:G\to\Diff_0^1([0,1])$ without fixed points in $(0,1)$, we have that the restriction of the action $\rho$ to $(0,1)$ is semi-conjugate to the standard action on $X$. 
\end{thm}

The condition of being locally moving means that for any open subinterval $I\subset X$, the subgroup $G_I$ of elements of $G$ fixing the complement $X\setminus I$ pointwise, acts on $I$ without fixed points. The condition of having independent groups of germs, in the context of the groups $F_\alpha$, means that $F_\alpha$ contains elements $f,g\in F_\alpha$ such that $Df(0)=\alpha$, $Dg(1)=\alpha$, and whose supports are disjoint. Both conditions are clearly satisfied in the groups $F_\alpha$. We finally deduce Corollary \ref{cor.Falpha}.

\appendix

\section{Direct computation}

In the proof of our main technical result (Theorem \ref{t.keythm}), we have used a criterion by Adouani and Marzougui in \cite{MR2427509}, for general piecewise $C^r$ homeomorphisms to be conjugate to an irrational rotation. As the general computation involves Faà di Bruno's formul\ae, it may be hard to digest for a non-specialist reader. In the case of low regularity ($r<4$), one can alternatively use the cocycle properties of the non-linearity and the Schwarzian derivative, to give a quite clean proof. 
We detail this in the case of low regularity ($r<3$).

\begin{lem}\label{l.reg_conjugation}
With notation as in the proof of Theorem \ref{t.keythm}, consider any diffeomorphism $\phi\in \Diff_0^\infty([\widetilde x,\widetilde z])$ satisfying the following properties: \begin{itemize}
	\item $\sigma(\phi)(\widetilde x)=D_+\phi(\widetilde x)/D_-\phi(\widetilde z)=\sigma(\widetilde T)(\widetilde x)$,
	\item $\phi(\widetilde y)=\widetilde y$.
\end{itemize}
Then $\phi \widetilde T \phi^{-1}$ is of class $C^{2}$.
Moreover, if the maps $\widetilde f$ and $\widetilde g$ are of class $C^{2+\tau}$, with $\tau\in (0,1]$, and $\phi$ satisfies
\begin{multline}\label{eq.conditionC2}
	D_-\widetilde T(\widetilde y)D_-\log D_-\phi(\widetilde x)-D_+\widetilde T(\widetilde y)D_+\log D_+\phi(\widetilde x)\\=D_-\widetilde T(\widetilde y)D_-\log D_-\widetilde T(\widetilde x)-D_+\widetilde T(\widetilde y)D_+\log D_+\widetilde T(\widetilde x),
\end{multline}
then $\phi \widetilde T \phi^{-1}$ is of class $C^{2+\tau}$. 
\end{lem}


\begin{proof}
	We first prove that $\phi \widetilde T \phi^{-1}$ is of class $C^1$.
	By the chain rule, we have
	\[
	D_{\pm}(\phi \widetilde T \phi^{-1})\circ \phi=\frac{D_\pm\phi\circ \widetilde T}{D_\pm \phi}D_\pm \widetilde T.
	\]
	As $D\widetilde T$ is continuous at every $t\notin \{\widetilde x,\widetilde y\}$, and $D\phi$ is continuous at every $t\neq\widetilde x$, we have that $D(\phi \widetilde T \phi^{-1})$ is continuous at every $t\notin \{\phi(\widetilde x),\phi(\widetilde y)\}=\{\widetilde x,\widetilde y\}$.
	For $t=\widetilde x$, we have
	\[
	\sigma(\phi \widetilde T \phi^{-1})(\widetilde x)=\sigma(\phi \widetilde T \phi^{-1})(\phi(\widetilde x))=\frac{\sigma(\phi)(\widetilde T(\widetilde x))}{\sigma( \phi)(\widetilde x)}\sigma(\widetilde T)(\widetilde x)=\frac{1}{\sigma( \phi)(\widetilde x)}\sigma(\widetilde T)(\widetilde x)=1.
	\]
	For $t=\widetilde y$, we have
	\begin{multline}\label{eq.jump}
	\sigma(\phi \widetilde T \phi^{-1})(\widetilde y)=\sigma(\phi \widetilde T \phi^{-1})(\phi(\widetilde y))=\frac{\sigma(\phi)(\widetilde T(\widetilde y))}{\sigma( \phi)(\widetilde y)}\sigma(\widetilde T)(\widetilde y)\\
	=\sigma(\phi)(\widetilde x)\cdot \sigma(\widetilde T)(\widetilde y)=\sigma(\widetilde T)(\widetilde x)\cdot \sigma(\widetilde T)(\widetilde y).
	\end{multline}
	Remark that condition \eqref{i.epsilon} for a local rotation guarantees that the jump of $D\widetilde T$ at $\widetilde y$ is the inverse of the jump of  $D\widetilde T$ at $\widetilde x\sim \widetilde z$.
		Indeed, we have 
		\[\sigma(\widetilde T)(\widetilde x)=D_+\widetilde f(\widetilde g(\widetilde y))/D_-\widetilde g(\widetilde f(\widetilde y)),\]
		whereas
		\[\sigma(\widetilde T)(\widetilde y)=D_+\widetilde g(\widetilde y)/D_-\widetilde f(\widetilde y),\]
		so that 
		\begin{equation}\label{eq.total_jump}
			\sigma(\widetilde T)(\widetilde x)\cdot \sigma(\widetilde T)(\widetilde y)=\frac{D_+\widetilde f(\widetilde g(\widetilde y))}{D_-\widetilde g(\widetilde f(\widetilde y))}\cdot \frac{D_+\widetilde g(\widetilde y)}{D_-\widetilde f(\widetilde y)}=\frac{D_+(\widetilde f\widetilde g)(\widetilde y)}{D_-(\widetilde g\widetilde f)(\widetilde y)}=1.
			\end{equation}
		Note that the last equality holds because of condition \eqref{i.epsilon} in the definition of a local rotation, which says that $fg=gf$ on a neighborhood of $y$, so after conjugation, the same holds for $\widetilde f$ and $\widetilde g$ on a neighborhood of $\widetilde y$.
	Combining \eqref{eq.total_jump} with \eqref{eq.jump}, we deduce that $\sigma(\phi \widetilde T \phi^{-1})(\widetilde y)=1$.
	This proves that $\phi \widetilde T \phi^{-1}$ is of class $C^1$. We next want to prove that $\phi \widetilde T \phi^{-1}$ is of class $C^2$, and for this we will show that the \emph{nonlinearity} $\cN(\phi \widetilde T \phi^{-1})$ is continuous. Recall that for a map $h$ of class $C^2$, the nonlinearity is defined by $\cN(h)=D\log Dh=\frac{D^2h}{Dh}$, and it satisfies the cocycle relation
	\begin{equation}\label{eq.cocycle}
	\cN(h_1h_2)=\cN(h_1)\circ h_2\cdot Dh_2+\cN(h_2).
	\end{equation}
	For the verifications below, we will write $\cN_\pm(h)=D_\pm\log D_\pm h$; it is straightforward to check that they still satisfy the cocycle relation \eqref{eq.cocycle}. Using this, we have
	\[
	\cN_\pm(\phi \widetilde T \phi^{-1})\circ \phi=\frac{1}{D_\pm\phi}\left (\cN_\pm(\phi)\circ \widetilde T\cdot D_\pm \widetilde T+\cN_\pm(\widetilde T)-\cN_\pm(\phi)\right ).
	\]
	As for the first derivative, we thus have that $\cN(\phi \widetilde T \phi^{-1})$ is continuous at every $t\notin \{\widetilde x,\widetilde y\}$. Continuity at the points $\widetilde x$ and $\widetilde y$ are equivalent, respectively, to the equations
	\begin{multline*}
	\sigma(\widetilde T)(\widetilde y)\left (
	\cN(\phi)(\widetilde T(\widetilde x))D_+\widetilde T(\widetilde x)+\cN_+(\widetilde T)(\widetilde x)-\cN_+(\phi)(\widetilde x)\right )\\=\cN(\phi)(\widetilde T(\widetilde x))D_-\widetilde T(\widetilde{x})+\cN_-(\widetilde{T})(\widetilde{x})-\cN_-(\phi)(\widetilde{x})
	\end{multline*}
	and
	\[
	\cN_+(\phi)(\widetilde x)D_+\widetilde T(\widetilde y)+\cN_+(\widetilde T)(\widetilde y)-\cN(\phi)(\widetilde y)=
	\cN_-(\phi)(\widetilde x)D_-\widetilde T(\widetilde{y})+\cN_-(\widetilde{T})(\widetilde{y})-\cN(\phi)(\widetilde{y}).
	\]
	Using that
	$\sigma(\widetilde T)(\widetilde {y})D_+\widetilde T(\widetilde x)=D_-\widetilde{T}(\widetilde{x})$ (because of \eqref{eq.total_jump}),
	the two equations are reduced to the linear system in the variables $N_-:=\cN_-(\phi)(\widetilde x)$ and $N_+:=\cN_+(\phi)(\widetilde x)$:
	\[
	\left\{\begin{array}{l}
		D_-\widetilde T(\widetilde y)N_--D_+\widetilde T(\widetilde y)N_+=D_-\widetilde T(\widetilde y)\cN_-(\widetilde{T})(\widetilde{x})-D_+\widetilde{T}(\widetilde{y})\cN_+(\widetilde T)(\widetilde{x}),\\[.5em]
		D_-\widetilde T(\widetilde y)N_--D_+\widetilde T(\widetilde y)N_+=-\cN_-(\widetilde T)(\widetilde{y})+\cN_+(\widetilde{T})(\widetilde{y}).
	\end{array} \right.
	\]
	Because of condition \eqref{eq.conditionC2}, the system has a solution if and only if
	\begin{equation*}\label{eq.conditionC2-1}
		D_-\widetilde T(\widetilde y)\cN_-(\widetilde{T})(\widetilde{x})-D_+\widetilde{T}(\widetilde{y})\cN_+(\widetilde T)(\widetilde{x})=-\cN_-(\widetilde T)(\widetilde{y})+\cN_+(\widetilde{T})(\widetilde{y}).
	\end{equation*}
	Using the cocycle relation \eqref{eq.cocycle}, this is equivalent to
	\begin{equation*}\label{eq.conditionC2-2}
		\cN_-(\widetilde{T}^2)(\widetilde{y})=\cN_+(\widetilde T^2)(\widetilde{y}).
	\end{equation*}
	In turns, this corresponds to the condition
	\begin{equation*}\label{eq.conditionC2-3}
		\cN_-(\widetilde{g}\widetilde{f})(\widetilde{y})=\cN_+(\widetilde f\widetilde g)(\widetilde{y}),
	\end{equation*}
	which holds because of condition \eqref{i.epsilon} in the definition of a local rotation, which says that $fg=gf$ on a neighborhood of $y$.
	
	The last statement is a direct consequence of the fact that a continuous function which is piecewise $\tau$-H\"older continuous is globally $\tau$-H\"older continuous.
\end{proof}

{\small \subsection*{Acknowledgments}
	We thank Joaqu\'in Brum, Nicol\'as Matte Bon, and Crist\'obal Rivas for discussions and interest around this project, and we thank the anonymous referee for the valuable comments.
	M.T. is partially supported by the project ANR
	Gromeov (ANR-19-CE40-0007), the project ANER
	Agroupes (AAP 2019 Région Bourgogne–Franche–Comté), and his host department IMB
	receives support from the EIPHI Graduate School (ANR-17-EURE-0002).}

\bibliographystyle{plain}

\bibliography{biblio.bib}

\begin{bibdiv}
\begin{biblist}

\bib{MR2427509}{article}{
      author={Adouani, Abdelhamid},
      author={Marzougui, Habib},
       title={Sur les hom\'{e}omorphismes du cercle de classe {$P\ C^r$} par
  morceaux {$(r\geq 1)$} qui sont conjugu\'{e}s {$C^r$} par morceaux aux
  rotations irrationnelles},
        date={2008},
        ISSN={0373-0956},
     journal={Ann. Inst. Fourier (Grenoble)},
      volume={58},
      number={3},
       pages={755\ndash 775},
         url={http://aif.cedram.org/item?id=AIF_2008__58_3_755_0},
      review={\MR{2427509}},
}

\bib{BieriStrebel}{book}{
      author={Bieri, Robert},
      author={Strebel, Ralph},
       title={On groups of {PL}-homeomorphisms of the real line},
      series={Mathematical Surveys and Monographs},
   publisher={American Mathematical Society, Providence, RI},
        date={2016},
      volume={215},
        ISBN={978-1-4704-2901-0},
         url={https://doi.org/10.1090/surv/215},
      review={\MR{3560537}},
}

\bib{CentralizersVn}{article}{
      author={{Bleak}, Collin},
      author={{Bowman}, Hannah},
      author={{Lynch}, Alison~Gordon},
      author={{Graham}, Garrett},
      author={{Hughes}, Jacob},
      author={{Matucci}, Francesco},
      author={{Sapir}, Eugenia},
       title={{Centralizers in the R. Thompson group \(V_n\)}},
        date={2013},
        ISSN={1661-7207},
     journal={{Groups Geom. Dyn.}},
      volume={7},
      number={4},
       pages={821\ndash 865},
}

\bib{BLT}{article}{
      author={Bonatti, Christian},
      author={Lodha, Yash},
      author={Triestino, Michele},
       title={Hyperbolicity as an obstruction to smoothability for
  one-dimensional actions},
        date={2019},
        ISSN={1465-3060},
     journal={Geom. Topol.},
      volume={23},
      number={4},
       pages={1841\ndash 1876},
         url={https://doi.org/10.2140/gt.2019.23.1841},
      review={\MR{3988090}},
}

\bib{BrinGuzman}{article}{
      author={Brin, Matthew~G.},
      author={Guzm\'{a}n, Fernando},
       title={Automorphisms of generalized {T}hompson groups},
        date={1998},
        ISSN={0021-8693},
     journal={J. Algebra},
      volume={203},
      number={1},
       pages={285\ndash 348},
         url={https://doi.org/10.1006/jabr.1997.7315},
      review={\MR{1620674}},
}

\bib{brum2021locally}{misc}{
      author={Brum, Joaquín},
      author={Bon, Nicolás~Matte},
      author={Rivas, Cristóbal},
      author={Triestino, Michele},
       title={Locally moving groups acting on the line and $\mathbb{R}$-focal
  actions},
        date={2021},
        note={arXiv:2104.14678},
}

\bib{Dzhalilov_2006}{article}{
      author={Dzhalilov, Akhtam},
      author={Liousse, Isabelle},
       title={Circle homeomorphisms with two break points},
        date={2006},
        ISSN={0951-7715},
     journal={Nonlinearity},
      volume={19},
      number={8},
       pages={1951\ndash 1968},
         url={https://doi.org/10.1088/0951-7715/19/8/010},
}

\bib{ghazouani2021priori}{misc}{
      author={Ghazouani, Selim},
      author={Ulcigrai, Corinna},
       title={{A priori bounds for GIETs, affine shadows and rigidity of
  foliations in genus 2}},
        date={2021},
        note={arXiv:2106.03529},
}

\bib{Ghys-bounded}{incollection}{
      author={Ghys, \'{E}tienne},
       title={Groupes d'hom\'{e}omorphismes du cercle et cohomologie
  born\'{e}e},
        date={1987},
   booktitle={The {L}efschetz centennial conference, {P}art {III} ({M}exico
  {C}ity, 1984)},
      series={Contemp. Math.},
      volume={58},
   publisher={Amer. Math. Soc., Providence, RI},
       pages={81\ndash 106},
         url={https://doi.org/10.24033/msmf.358},
      review={\MR{893858}},
}

\bib{GhysSergiescu}{article}{
      author={Ghys, \'{E}tienne},
      author={Sergiescu, Vlad},
       title={Sur un groupe remarquable de diff\'{e}omorphismes du cercle},
        date={1987},
        ISSN={0010-2571},
     journal={Comment. Math. Helv.},
      volume={62},
      number={2},
       pages={185\ndash 239},
         url={https://doi.org/10.1007/BF02564445},
      review={\MR{896095}},
}

\bib{Herman}{article}{
      author={Herman, Michael-Robert},
       title={Sur la conjugaison diff\'{e}rentiable des diff\'{e}omorphismes du
  cercle \`a des rotations},
        date={1979},
        ISSN={0073-8301},
     journal={Inst. Hautes \'{E}tudes Sci. Publ. Math.},
      number={49},
       pages={5\ndash 233},
         url={http://www.numdam.org/item?id=PMIHES_1979__49__5_0},
      review={\MR{538680}},
}

\bib{HmiliLiousse}{article}{
      author={{Hmili}, Hadda},
      author={{Liousse}, Isabelle},
       title={{Dynamique des \'echanges d'intervalles des groupes de
  Higman-Thompson \(V_{r, m}\)}},
    language={French},
        date={2014},
        ISSN={0373-0956},
     journal={{Ann. Inst. Fourier}},
      volume={64},
      number={4},
       pages={1477\ndash 1491},
}

\bib{hyde2021subgroups}{article}{
      author={Hyde, James},
      author={Moore, Justin~Tatch},
       title={{Subgroups of $\mathrm{PL}_+ I$ which do not embed into
  Thompson's group $F$}},
        date={to appear},
     journal={Groups Geom. Dyn.},
}

\bib{KhaninKocic}{incollection}{
      author={{Khanin}, Konstantin},
      author={{Koci\'c}, Sa\v{s}a},
       title={{Rigidity for a class of generalized interval exchange
  transformations}},
        date={2017},
   booktitle={{Dynamical systems, ergodic theory, and probability. In memory of
  Kolya Chernov. Conference dedicated to the memory of Nikolai Chernov,
  University of Alabama at Birmingham, Birmingham, AL, USA, May 18--20, 2015}},
   publisher={Providence, RI: American Mathematical Society (AMS)},
       pages={161\ndash 167},
}

\bib{KhaninTeplisnky}{article}{
      author={{Khanin}, Konstantin},
      author={{Teplinsky}, Alexey},
       title={{Herman's theory revisited}},
        date={2009},
        ISSN={0020-9910},
     journal={{Invent. Math.}},
      volume={178},
      number={2},
       pages={333\ndash 344},
}

\bib{KKL}{article}{
      author={Kim, Sang-hyun},
      author={Koberda, Thomas},
      author={Lodha, Yash},
       title={Chain groups of homeomorphisms of the interval},
        date={2019},
        ISSN={0012-9593},
     journal={Ann. Sci. \'{E}c. Norm. Sup\'{e}r. (4)},
      volume={52},
      number={4},
       pages={797\ndash 820},
         url={https://doi.org/10.24033/asens.239},
      review={\MR{4038452}},
}

\bib{liousse}{article}{
      author={{Liousse}, Isabelle},
       title={{Nombre de rotation, mesures invariantes et ratio set des
  hom\'eomorphismes affines par morceaux du cercle}},
    language={French},
        date={2005},
        ISSN={0373-0956},
     journal={{Ann. Inst. Fourier}},
      volume={55},
      number={2},
       pages={431\ndash 482},
}

\bib{Coherent}{article}{
      author={Lodha, Yash},
       title={Coherent actions by homeomorphisms on the real line or an
  interval},
        date={2020},
        ISSN={0021-2172},
     journal={Israel J. Math.},
      volume={235},
      number={1},
       pages={183\ndash 212},
         url={https://doi.org/10.1007/s11856-019-1954-7},
      review={\MR{4068782}},
}

\end{biblist}
\end{bibdiv}

\medskip

\noindent\textit{Michele Triestino\\
	Institut de Math\'ematiques de Bourgogne (IMB, UMR CNRS 5584)\\
	Universit\'e Bourgogne Franche-Comt\'e\\
	9 av.~Alain Savary, 21000 Dijon, France\\}
\href{mailto:michele.triestino@u-bourgogne.fr}{michele.triestino@u-bourgogne.fr}


\end{document}